\documentclass[centertags,12pt]{amsart}
\usepackage{latexsym}
\usepackage{amssymb}
\newtheorem{claim}{Claim}
\newtheorem*{claim*}{Claim}
\newtheorem{theorem}{Theorem}
\newtheorem{lemma}{Lemma}
\theoremstyle{remark}
\newtheorem*{remark*}{Remark}
\newtheorem*{remarks*}{Remarks}
\def\N{\mathbb N}
\def\P{\mathbb P}
\def\cD{\mathcal D}
\def\cJ{\mathcal J}
\def\cO{\mathcal O}
\def\fq{\mathbb F_q}

\def\df{{\rm div}}
\def\dinf{{\rm div}_{\infty}}
\def\deg{{\rm deg}}
\def\sq{\sqrt{q}}
\def\frj{{\rm Fr}_{{\mathcal J}}}
\def\frx{{\rm Fr}_{X}}
\begin{document}
\author[F.~Torres]{Fernando Torres}\thanks{The paper was partially written
while the author was visiting the University of Essen supported by a grant
from the Graduiertenkolleg ``Theoretische und experimentelle Methoden der
Reinen Mathematik". The results of this paper were announced in
\cite{t}}
\title[A Deligne-Lusztig curve]{The Deligne-Lusztig curve associated\\
to the Suzuki group} 
\address{ICTP, Mathematics Section,  P.O. Box 586, 34100, Trieste, Italy}
\email{ftorres@ictp.trieste.it}
\begin{abstract}
We give a characterization of the Deligne-Lusztig curve  
associated to the Suzuki group $Sz(q)$  
based on the genus and the number of $\fq$-rational points of the curve. 
\end{abstract}
\maketitle
\noindent {\bf \S0.} Throughout this paper by a curve we mean a projective,
geometrically irreducible, and non-singular algebraic curve defined over
the finite field $\fq$ with $q$ elements.

Let $N_q(g)$ denote the maximum number of $\fq$-rational points that 
a curve of genus $g$ can have. The number $N_q(g)$ is bounded from above
by the Hasse-Weil bound. A curve $X$ of genus $g$ is said to be {\it
optimal} over $\fq$, if its number of $\fq$-rational points $\# X(\fq)$ is
equal to $N_q(g)$. Optimal curves are very useful for applications to
coding theory via Goppa construction \cite{go}.

Besides {\it maximal curves} over $\fq$, that is, curves whose number of 
$\fq$-rational points attains the Hasse-Weil bound and some curves of
small genus, the only known examples of optimal curves are the 
Deligne-Lusztig curves associated to the Suzuki group $Sz(q)$ and to the
Ree 
group $R(q)$ \cite[\S11]{dl}, \cite{h}. Arithmetical and
geometrical 
properties of maximal curves were studied in 
\cite{rsti}, \cite{i}, \cite{gvl} (see the references therein), \cite{ft}
and \cite{fgt}. 
The other mentioned optimal curves were studied in 
\cite{hsti}, \cite{h}, \cite{p} and \cite{hp}. Hansen and Pedersen 
\cite[Thm. 1]{hp} stated the uniqueness, up to $\fq$-isomorphism, of the curve
corresponding to $R(q)$ based on the genus, the number of
$\fq$-rational points, and the group of $\fq$-automorphisms of the curve. 
They observed a similar result for the curve corresponding to $Sz(q)$   
(cf. \cite[p.100]{hp}) as a consequence of its uniqueness up to $\bar
\fq$-isomorphism (cf. \cite{he}). Hence, by \cite{hsti}, the curve under
study in this paper is $\fq$-isomorphic to the plane curve given by 
$$
y^q-y=x^{q_0}(x^q-x),
$$
where $q_0=2^s$ and $q=2q_0$. Its genus is $q_0(q-1)$ and its number of
$\fq$-rational points is $q^2+1$ (loc. cit.).

In this note, we will be mainly concerned with the uniqueness, up to 
$\fq$-isomorphism, of  the Deligne-Lusztig curve associated to $Sz(q)$   
based on the genus and the number of 
$\fq$-rational points of the curve only. Our main result is
the following 
\begin{theorem}\label{thm1} For $s\in \mathbb N$, set $\ q_0:= 2^{2s}$ and  
$q:= 2q_0$. Let $X$ be a curve over $\fq$ of genus $g$ satisfying the 
following hypotheses:
\begin{enumerate}
\item $g=q_0(q-1)$;
\item $\#X(\fq)=q^2+1$.
\end{enumerate} 
Then $X$ is $\fq$-isomorphic to the Deligne-Lusztig curve associated to
the Suzuki group $Sz(q)$.
\end{theorem}
We also state an analogous result for the Hermitian curve over $\fq$, 
that is,   
the curve defined by $y^{\sq}+y=x^{\sq+1}$, with $q$ a square. This  is a 
maximal curves over $\fq$ whose genus is  the biggest genus 
that a maximal curve over $\fq$ can have, namely $\sq(\sq-1)/2$ (cf. 
\cite{ih}). Several 
characterizations of Hermitian curves have already been given, see e.g. 
\cite{hstv} (and the references therein), \cite{rsti} and \cite{ft}. We 
remark that a Hermitian curve is a Deligne-Lusztig curve associated
to a 
projective special linear group (cf. \cite{h}). It is unique up to 
$\fq$ (cf. \cite[p.100]{hp}) because of its 
uniqueness up to $\bar \fq$ (cf. \cite{sti}). We also remark that Theorem 
\ref{thm2} below follows from the results in \cite{ft} (here, we give a 
new proof).  
\begin{theorem}\label{thm2}
Let $X$ be a curve over $\fq$ of genus $g$ satisfying the
following hypotheses:
\begin{enumerate}
\item $g>(\sq-1)^2/4$;
\item $\#X(\fq)=q+2g\sq+1$ (i.e. $X$ is maximal over $\fq$).
\end{enumerate} 
Then $X$ is $\fq$-isomorphic to the Hermitian curve over $\fq$. 
In particular\\
$g=\sq(\sq-1)/2$.
\end{theorem}
To prove these theorems we use geometrical properties of the Zeta Function
of curves \cite[\S1]{fgt} and a particular case of the theory of 
Frobenius orders  
associated to linear system of curves which was developed by 
St\"ohr and Voloch \cite{sv} for an improvement of the Hasse-Weil bound. 

It is a pleasure to thank: R. Pellikaan and H. Stichtenoth for useful 
conversations; A. Cossidente for having let me include his observation 
in this paper (see Appendix). There he points out a relation between the 
Deligne-Lusztig curve associated to the Suzuki group and the Suzuki-Tits
ovoid. I also thank Prof. J.F. Voloch for his interest in this work.   
\medskip

\noindent {\bf \S1. A lemma on non-gaps.} Let $X$ be a curve satisfying 
the
hypotheses of Theorem \ref{thm1}, and $\frj$ the Frobenius endomorphism 
(relative to $\fq$) of the Jacobian $\cJ={\cJ}_X$ of $X$. Then, by 
Serre-Weil's explicit formulae 
(cf. \cite{se}, \cite{h}), 
$$ 
h(t)=(t^2+2q_0t+q)^g 
$$ 
is the 
characteristic polynomial of $\frj$. Furthermore, 
$$ 
\frj^2+2q_0\frj+q^2I=0\quad \mbox{on}\ \
\cJ,  
$$ 
due to the 
semi-simplicity of $\frj$ and the fact that the representation of 
endomorphisms of $\cJ$ on the Tate module is faithful (cf. \cite[Thm. 
2]{ta}, \cite[VI, \S3]{l}). 
 
Let $P_0\in X(\fq)$ and consider the natural map $f=f^{P_0}:X\to
\cJ$ giving by \\ 
$P\mapsto [P-P_0]$. By using $f\circ\frx=\frj\circ f$, 
where $\frx$ stands for the Frobenius morphism on $X$ relative to $\fq$, 
we then have 
\begin{equation}\label{eq1} 
\frx^2(P)+2q_0\frx(P)+qP\sim(q+2q_0+1)P_0,\quad \mbox{for each}\ \ P\in X.  
\end{equation} 
It turns out that $X$ is equipped with the linear system 
$$
\cD=g^r_{q+2q_0+1}:=|(q+2q_0+1)P_0|, 
$$ 
which is independent of $P_0\in
X(\fq)$ by (\ref{eq1}).  
\begin{lemma}\label{lemma1} Let $X$ be a curve satisfying the hypotheses 
of Theorem \ref{thm1}. For each $Q\in X(\fq)$, the first positive non-gap  
at $Q$ is $q$.
\end{lemma}
\begin{proof} Let $x: X\to \mathbb P^1(\bar \fq)$ be such that $\dinf
(x)=mQ$, $m$  
being the first positive non-gap at $Q$. We show firstly that $q$
is a non-gap at a generic point of 
$X$ (\cite[proof of Prop. 1.8(i)]{fgt}). Applying (as in \cite[IV, Ex. 
2.6]{har}) $\frx$ to the equivalence in (\ref{eq1}) we obtain 
$$
\frx^3(P)+(2q_0-1)\frx^2(P)+(q-2q_0)\frx(P)\sim qP.  
$$ 
Consequently, if $\frx^i(P)\neq P$, $i=1,2,3$, $q$ is a non-gap at
$P$. In particular 
$$ 
m\le m_1(P)\le q, 
$$ 
where $P$ is a non-Weierstrass point and $m_1(P)$ is the first positive 
non-gap at $P$. On the other hand, $\# X(\fq)\le 1+qm$ (see e.g. 
\cite[Thm.  1(b)]{le}) 
and so $m\ge q$. Now the result follows.  
\end{proof} 

{\bf \S2. The linear system $\cD$.} Let $X$ be a curve of genus $g$ 
satisfying the hypothesis of the theorem. We are going to apply \cite{sv} 
to the linear system $\cD$ introduced in \S1. The key property of $\cD$
will be the fact that $\frx(P)$ belongs to the tangent line at $P$ for a generic  
$P$ (see Claim \ref{claim3}(1)).

For $P\in X$, let $0=j_0(P)<j_1(P)<\ldots<j_r(P)$ be the $(\cD,P)$-orders 
(\cite[\S1]{sv}) and $m_i(P)$ the $i$th non-gap at $P$ 
($m_0(P):=0$). If $P\in X(\fq)$, from (\ref{eq1}) we have 
$m_r(P)=q+2q_0+1$ and 
\begin{equation}\label{eq2}
j_i(P)=m_r(P)-m_{r-i}(P)\quad \mbox{for}\ \ i=0,\ldots,r.  
\end{equation}
Then by Lemma \ref{lemma1} and (\ref{eq2}), 
\begin{equation}\label{eq3}
j_r=(P)=q+2q_0+1,\quad j_{r-1}(P)=2q_0+1\qquad \mbox{for}\ \  P\in X(\fq).
\end{equation} 
Let $\epsilon_0=0<\epsilon_1<\ldots<\epsilon_r$ be the 
orders of $\cD$ (i.e. the $(\cD,P)$-orders for $P$ generic). From
(\ref{eq1}) we see that $1, 2q_0, q$ are orders for $\cD$. Thus $r\ge 3$ 
and $\epsilon_1=1$.  Therefore
\begin{equation}\label{eq4}
\epsilon_r=q\qquad\text{ and}\qquad 2q_0\le\epsilon_{r-1}\le 1+2q_0, 
\end{equation} 
because $\epsilon_i\le j_i(P)$ (\cite[p.5]{sv}).   
Let $\nu_0=0<\nu_1<\ldots<\nu_{r-1}$ be the 
$\fq$-Frobenius orders of $\cD$ (which is obtained from the order sequence
by deleting an specific order; cf. \cite[Prop. 2.1]{sv}). These numbers 
satisfy \cite[Corollary 2.6]{sv} 
\begin{equation}\label{eq5} 
\nu_i\le j_{i+1}(P)-j_1(P)\quad \text{whenever}\ \ P\in X(\fq).
\end{equation}
From (\ref{eq1}) we have that $\frx(P)$ belongs to the osculating 
hyperplane at $P$ provided that $P\not\in X(\fq)$. Consequently 
$\nu_{r-1}=\epsilon_r=q$.  
\begin{claim}\label{claim1}\quad $\epsilon_{r-1}=2q_0$.  
\end{claim}
\begin{proof} Suppose that $\epsilon_{r-1}> 
2q_0$. Then, by (\ref{eq4}), $\epsilon_{r-2}=2q_0$ and 
$\epsilon_{r-1}=2q_0+1$. Moreover, by (\ref{eq5}) and (\ref{eq3}), 
$\nu_{r-2}\le 2q_0=\epsilon_{r-2}$. Hence the $\fq$-Frobenius orders of 
$\cD$ would be $\epsilon_0,\epsilon_1,\ldots,\epsilon_{r-2}$,  
and $\epsilon_r=q$. 
 
Associated to $\cD$ there exists an effective divisor $S$ on $X$ such that  
$$ 
{\rm deg}(S)=\sum_{i=0}^{r-1}\nu_i (2g-2)+(q+r)(q+2q_0+1).
$$ 
For $P\in X(\fq)$, we have 
\begin{equation}\label{eq6}
v_P(S)\ge \sum_{i=1}^{r}(j_i(P)-\nu_{i-1})\ge
(r-1)j_1(P)+1+2q_0\ge r+2q_0,
\end{equation} 
as follows from  \cite[Prop. 2.4, Corollary 2.6]{sv}, (\ref{eq3}) and $\nu_{r-1}=q$. 
Thus $\deg(S)\ge (r+2q_0)\#X(\fq)$, and
from the formula for $\deg(S)$, $2g-2=(2q_0-2)(q+2q_0+1)$, and 
$\#X(\fq)=q^2+1=(q-2q_0+1)(q+2q_0+1)$ it follows that 
$$
\sum_{i=1}^{r-2}\nu_i=\sum_{i=1}^{r-2}\epsilon_i\ge (r-1)q_0.  
$$ 
Now, as $\epsilon_i+\epsilon_j\le \epsilon_{i+j}$ for $i+j\le r$ 
\cite[Thm. 1]{e}, then we would have 
$$ 
(r-1)\epsilon_{r-2}\ge 2\sum_{i=0}^{r-2}\epsilon_i,  
$$ 
and hence $\epsilon_i+\epsilon_{r-2-i}=\epsilon_{r-2}$ for
$i=0,\ldots,r-2$. In particular,
$\epsilon_{r-3}=2q_0-1$ and by the $p$-adic criterion (cf. \cite[Corollary
1.9]{sv} we would have $\epsilon_i=i$ for $i=0,1,\ldots,r-3$. These facts
imply $r=2q_0+2$. 

Finally, we are going to see that this is a contradiction via 
Castelnuovo's genus bound for curves in projective spaces. We 
notice that $\nu_{r-2}=2q_0$, (\ref{eq3}) and (\ref{eq5}) imply  
$j_1(P)=1$ for $P\in X(\fq)$. Then, by (\ref{eq2}),  
$m_{r-1}(P)=q+2q_0$ 
and hence $\cD$ is simple. Therefore Castelnuovo's formula (see \cite{c},
\cite[p. 116]{acgh}, \cite[Corollary 2.8]{ra}) applied to $\cD$ implies
$$
2g=2q_0(q-1)\le \frac{(q+2q_0-(r-1)/2)^2}{r-1}.
$$
For $r=2q_0+2$ this gives $2q_0(q-1)< (q+q_0)^2/2q_0=q_0q+q/2+q_0/2$, a
contradiction.
\end{proof}
\begin{claim}\label{claim2}
There exists $P_1\in X(\fq)$ such that 
$$
\left\{ \begin{array}{ll}
j_1(P_1)=1         & {} \\
j_i(P_1)=\nu_{i-1}+1 & \mbox{if}\ i=2,\ldots, r-1.
\end{array}\right.
$$
\end{claim}
\begin{proof} By (\ref{eq6}), it is enough to show that there exists
$P_1\in X(\fq)$ 
such that $v_{P_1}(S)=r+2q_0$. Suppose that $v_P(S)\ge r+2q_0+1$ for each 
$P\in X(\fq)$. Then from the formula for $\deg(S)$ we would have that 
$$
\sum_{i=0}^{r-1}\nu_i \ge q+rq_0+1,
$$
and, as $\epsilon_1=1$, $\nu_{r-1}=q$ and $\nu_i\le \epsilon_{i+1}$, that 
$$
\sum_{i=0}^{r-1}\epsilon_i \ge rq_0+2.
$$
By \cite[Thm. 1]{e}, we then conclude that 
$r\epsilon_{r-1}\ge 2rq_0+4$, i.e. $\epsilon_{r-1}>2rq_0$, a contradiction
(cf. Claim \ref{claim1}).
\end{proof}
\begin{claim}\label{claim3}
\begin{enumerate} 
\item \ $\nu_1>\epsilon_1=1$.
\item \ $\epsilon_2$ is a power of two.
\end{enumerate}
\end{claim}
\begin{proof} Statement (2) is consequence of the $p$-adic criterion 
\cite[Corollary 1.9]{sv}.  
Suppose that $\nu_1=1$. Then by Lemma \ref{lemma1}, Claim \ref{claim2}, (\ref{eq3})
and  
(\ref{eq2}) there would be a point $P_1\in X$ such that the Weierstrass
semigroup at  
$P_1$ would contain the semigroup $H:=\langle q, q+2q_0-1, q+2q_0, 
q+2q_0+1\rangle $. Then $g\le  \#(\N\setminus
H)$, a contradiction (cf. Remark below). 
\end{proof}
\begin{remark*} Let $H$ be the semigroup defined above. We are going to
show that $\tilde g:= \#(\N \setminus H)=g-q_0^2/4$. 
To start with, notice that  
$L:=\cup_{i=1}^{2q_0-1}L_i$ is a complete system of residues
module $q$, where
$$
\begin{array}{lll}
L_i & = &
\{iq+i(2q_0-1)+j: j=0,\ldots,2i\}\quad  \mbox{if}\ \ 1\le i\le q_0-1,\\
L_{q_0} & = & \{q_0q+q-q_0+j:j=0,\ldots,q_0-1\},\\
L_{q_0+1} & = & \{(q_0+1)q+1+j:j=0,\ldots,q_0-1\},\\ 
L_{q_0+i} & = & 
\{(q_0+i)q+(2i-3)q_0+i-1+j: 
j=0,\ldots,q_0-2i+1\}\cup\\
          &   & \{(q_0+i)q+(2i-2)q_0+i+j: j=0,\ldots q_0-1\}\quad  
\mbox {if}\ \ 2\le i\le q_0/2,\\
L_{3q_0/2+i} & = & 
\{(3q_0/2+i)q+(q_0/2+i-1)(2q_0-1)+q_0+2i-1+j:\\
             &  &  
j=0,\ldots,q_0-2i-1\}\quad \mbox {if}\ \ 1\le i\le q_0/2-1.
\end{array}
$$
Moreover, for each $\l \in L$, $\l \in H$ and $\l-q\not\in H$. Hence
$\tilde g$ can be computed by summing up the coefficients of $q$ from the
above list (see e.g. \cite[Thm. p.3]{sel}), i.e.
$$
\begin{array}{lll}
\tilde g & = & \sum_{i=1}^{q_0-1}i(2i+1)+q_0^2+(q_0+1)q_0+
\sum_{i=2}^{q_0/2}(q_0+i)(2q_0-2i+2)+\\
         &   & 
\sum_{i=1}^{q_0/2-1}(3q_0/2+i)(q_0-2i)=q_0(q-1)-q_0^2/4.
\end{array}
$$
\end{remark*}
In the remaining part of this section let $P_1\in X(\fq)$ be a point  
satisfying Claim \ref{claim2}. We set $m_i:= m_i(P_1)$ and denote by $v$   
the valuation at $P_1$. 

By Claim \ref{claim3}(1) we have $\nu_i=\epsilon_{i+1}$ for 
$i=1,\ldots, r-1$. Therefore from (\ref{eq2}),  
(\ref{eq3}) and Claim \ref{claim2},
\begin{equation}\label{eq7}
\left\{\begin{array}{ll}
m_i=2q_0+q-\epsilon_{r-i} & \mbox{if}\ i=1,\ldots r-2\\
m_{r-1}=2q_0+q,\ \ m_r=1+2q_0+q. & {}
\end{array}\right.
\end{equation}
Let $x, y_2,\ldots, y_r\in \fq(X)$ be such that $\dinf(x)=m_1P_1$, and 
$\dinf (y_i)=m_i P_1$ for $i=2,\ldots, r$. Let denote by $D^{(i)}$ the 
$i$th Hasse derivative defined by the separating variable $x$.   
The fact that $\nu_1>1$ means that the following 
matrix has rank two (cf. \cite[proof of Prop. 2.1]{sv})
$$
\left( \begin{array}{ccccc}
1 & x^q & y_2^q &\ldots &y_r^q\\
1 & x   & y_2   &\ldots &y_r\\
0 & 1   & D^{(1)}y_2   &\ldots& D^{(1)}y_r
\end{array} \right).
$$
In particular, 
\begin{equation}\label{eq8}
y_i^q-y_i= D^{(1)}y_i(x^q-x) \quad \text{for}\ \ i=2,\ldots, r.
\end{equation}
\begin{claim}\label{claim4} 
\begin{enumerate}
\item For $P\in X(\fq)$, the divisor $(2g-2)P$ 
is canonical, i.e. the Weierstrass semigroup at $P$ is symmetric. 
\item  Let $m$ be a non-gap at $P_1\in X(\fq)$. If $m<q+2q_0$, then $m\le
q+q_0$. 
\item For $i=2,\ldots,r$ there exists $g_i\in \fq(X)$ such that $
D^{(1)}y_i=g_i^{\epsilon_2}$. 
Furthermore, $\dinf(g_i)=\frac{qm_i-q^2}{\epsilon_2}P_1$.
\end{enumerate}
\end{claim}
\begin{proof} (1) Since $2g-2=(2q_0-2)(q+2q_0+1)$, by (\ref{eq1}) 
we can assume $P=P_1$.  Now the case $i=r$ of 
Eqs. (\ref{eq8}) implies $v(dx)=2g-2$ and we are done.

(2) The numbers $q, q+2q_0$ and $q+2q_0+1$ are non-gaps at $P_1$ (cf. 
(\ref{eq7})). Then the numbers
$$
(2q_0-2)q+q-4q_0+j\qquad j=0,\ldots,q_0-2
$$
are also non-gaps at $P_1$. Therefore from the symmetry of the
Weierstrass semigroup at $P_1$,
$$
q+q_0+1+j\qquad j=0,\ldots,q_0-2
$$
are gaps at $P_1$. Now the proof follows.

(3) Set $f_i:= D^{(1)}y_i$. By the product rule applied to (\ref{eq8}),\\  
$D^{(j)}y_i=(x^q-x)D^{(j)}f_i+D^{(j-1)}f_i$ for $1\le j<q$. 
Then, $D^{(j)}f_i=0$ for $1\le 
j<\epsilon_2$, because the matrices
$$
\left( \begin{array}{ccccc}
1 & x   & y_2   &\ldots &y_r\\
0 & 1   & D^{(1)}y_2   &\ldots& D^{(1)}y_r\\
0 & 0   & D^{(j)}y_2   &\ldots& D^{(j)}y_r
\end{array} \right),
\quad 2\le j<\epsilon_2
$$
have rank two.  
Therefore, as $\epsilon_2$ is a power of 2 (cf. Claim 
\ref{claim3}(2)), by \cite[Satz 10]{hasse}, 
$f_i=g_i^{\epsilon_2}$ for some $g_i\in \fq(X)$. 

Finally, from the proof of item (1) we have that $x-x(P)$ is a local 
parameter at $P$ if $P\neq P_1$. Then, by the election of the $y_i$'s, 
$g_i$ has no pole but in $P_1$, and from (\ref{eq8}),   
$v(g_i)=-(qm_i-q^2)/\epsilon_2$. 
\end{proof}
\begin{claim}\label{claim5}\quad $r=4$ and $\epsilon_2=q_0$.
\end{claim}
\begin{proof} 
We know  that $r\ge 3$. We claim that $r\ge 4$; in 
fact, if $r=3$ we would have $\epsilon_2=2q_0$, $m_1=q$, $m_2=q+2q_0$,
$m_3=q+2q_0+1$, and hence $v(g_2)=-q$ ($g_2$ being as in Claim 
\ref{claim4}(3)). Therefore, 
after some $\fq$-linear transformations, the case $i=2$ of 
(\ref{eq8}) reads 
$$
y_2^q-y_2=x^{2q_0}(x^q-x).
$$
Now the function $z:= y_2^{q_0}-x^{q_0+1}$ satisfies
$z^q-z=x^{q_0}(x^q-x)$ and we find that $q_0+q$ is
a non-gap at $P_1$ (cf. \cite[Lemma 1.8]{hsti}). This contradiction
eliminates the case $r=3$.

Let $r\ge 4$ and $2\le i\le r$. By Claim \ref{claim4}(3), 
$(qm_i-q^2)/\epsilon_2$ is a 
non-gap at $P_1$, and since $(qm_i-q^2)/\epsilon_2\ge m_{i-1}\ge q$, from 
(\ref{eq7}) we have 
$$
2q_0\ge \epsilon_2 +\epsilon_{r-i}\qquad \mbox{for}\ i=2,\ldots,r-2.
$$
In particular, $\epsilon_2\le q_0$. On the other hand, by Claim 
\ref{claim4}(2) we must have $m_{r-2}\le q+q_0$ and so, by
(\ref{eq7}), we find that $\epsilon_2\ge q_0$, i.e. $\epsilon_2=q_0$. 

Finally we show that $r=4$. $\epsilon_2=q_0$ 
implies $\epsilon_{r-2}\le q_0$. Since $m_2\le q+q_0$ (cf. Claim
\ref{claim4}(2)), by (\ref{eq7}), we have $\epsilon_{r-2}\ge q_0$.
Therefore $\epsilon_{r-2}=q_0=\epsilon_2$, i.e. $r=4$.
\end{proof}
{\bf \S3. Proof of Theorem \ref{thm1}.} Let $P_1\in X(\fq)$ be as above. By
(\ref{eq8}), Claim \ref{claim4}(3) and 
Claim \ref{claim5} we have the following equation   
$$
y_2^q-y_2=g_2^{q_0}(x^q-x),
$$
where $g_2$ has no pole except at $P_1$. Moreover from (\ref{eq7}) 
we have that   
$m_2=q_0+q$ and hence that $v(g_2)=-q$ (cf. \ref{claim4}(3)). Thus 
$g_2=ax+b$ with $a,b\in \fq$, $a\neq 0$, and after some $\fq$-linear
transformations we obtain Theorem \ref{thm1}.
\begin{remarks*} (1) From the above computations we conclude that the
Deligne-Lusztig 
curve associated to the Suzuki group (as in \S0) is equipped with a complete and
simple base-point-free
linear system of dimension 4 and degree $q+2q_0+1$, namely
$\cD=|(q+2q_0+1)P_0|$, $P_0\in X(\fq)$. Such a linear system is an
$\fq$-invariant. The orders of $\cD$ (resp. the $\fq$-Frobenius orders)
are $0, 1, q_0, 2q_0$ and $q$ (resp. $0, q_0, 2q_0$ and $q$).

(2) There exists $P_1\in X(\fq)$ such that the $(\cD,P_1)$-orders are
$0,1,q_0+1,
2q_0+1$ and $q+2q_0+1$ (cf. Claim \ref{claim2}, (\ref{eq2}) and 
(\ref{eq3})). Now 
we show that the above sequence is, in fact, the
$(\cD,P)$-orders for each $P\in X(\fq)$. To see this, notice that
$$
\deg(S)=(3q_0+q)(2g-2)+(q+4)(q+2q_0+1)=(4+2q_0)\#X(\fq).
$$
Let $P\in X(\fq)$. From (\ref{eq6}), we conclude that 
$v_P(S)=\sum_{i=1}^{4}(j_i(P)-\nu_{i-1})=4+2q_0$ and hence, by (\ref{eq5}),
that $j_1(P)=1$, $j_2(P)=q_0+1$, $j_3(P)=2q_0+1$, and $j_4(P)=q+2q_0+1$.

(3) Then, by (\ref{eq2}), the Weierstrass semigroup $H(P)$ at $P\in
X(\fq)$   
contains the 
semigroup $H:= \langle q,q+q_0,q+2q_0,q+2q_0+1\rangle$. Indeed,  
$H(P)=H$ since $\#(\mathbb N\setminus H)=g=q_0(q-1)$ (this can be 
proved as in the remark after Claim \ref{claim3}; see also 
\cite[Appendix]{hsti}). 

(4) Let $R$ be the ramification divisor of $\cD$ \cite[\S1]{sv}. Its
support is the set of $\cD$-Weierstrass points of $X$ and  
$$
\deg(R)=\sum_{i=0}^{4}\epsilon_i(2g-2)+5(q+2q_0+1).
$$
If $P\in X(\fq)$, $v_P(R)=2q_0+3$, as follows from items (1), (2) and
\cite[Thm. 1.5]{sv}. We find that $\deg(R)=(2q_0+3)\#X(\fq)$ and 
consequently the set of $\cD$-Weierstrass points of $X$ is equal to 
$X(\fq)$. In 
particular, the $(\cD,P)$-orders for $P\not\in X(\fq)$ are $0, 1, q_0,  
2q_0$ and $q$. 

(5) We can use the above computations to obtain information on orders
for the canonical morphism. By using the fact that $(2q_0-2)\cD$ is
canonical (cf. Claim 
\ref{claim4}(1)) and item (4), we see that the set 
$$
\{a+q_0b+2q_0c+qd: a+b+c+d \le 2q_0-2\}
$$
is contained in the set of orders (for the
canonical
morphism) at non-rational points. (By considering certain rational
functions of 
the curve, similar computations were 
obtained in \cite[\S4]{gsti}.)

(6) Finally, we remark that $X$ is 
non-classical for the canonical morphism. We state two proofs of this fact  
(cf. loc. cit. and \cite[Prop.
1.8]{fgt}): by item (5), the number $q_0q$ is an order and since
$q_0q>g-1$ the curve is non-classical. On the other hand, from the proof
of Lemma \ref{lemma1}, $q$ is a non-gap at a generic point $P$. Since
$g\ge q$, this also implies the non-classicality of $X$.
\end{remarks*}
{\bf \S4. Proof of Theorem \ref{thm2}.} The proof will be similar to the
proof of Theorem \ref{thm1}. Let $X$ be a curve satisfying the
hypotheses of Theorem \ref{thm2}. 

It is well known that the characteristic  
polynomial of the Frobenius endomorphism of the Jacobian of a maximal curve over
$\fq$ of genus $g$ is 
$h(t)=(t+\sq)^{2g}$. Then, as in \S1, $X$ is equipped with the linear system    
$\cD:=g^r_{\sq+1}=|(\sq+1)P_0|$ with $P_0\in X(\fq)$. Here for each $P\in
X$, 
\begin{equation}\label{last}
\sq P+\frx(P)\sim (\sq+1)P_0.
\end{equation}
As in \S1, let $j_0(P)=0<\ldots<j_r(P)$ be the 
$(\cD,P)$-orders at $P$, $\epsilon_0=0,\ldots,\epsilon_r$ (resp. 
$\nu_0=0,\ldots,\nu_{r-1}$)  
the orders of $\cD$ (resp. the $\fq$-Frobenius orders of $\cD$). By 
(\ref{last}), $\nu_{r-1}=\epsilon_r$, $j_r(P)=m_r(P)=\sq+1$ for each $P\in 
X(\fq)$, and 1 and $\sq$ are 
$(\cD,P)$-orders for each $P\not\in X(\fq)$. Then $j_1(P)=1$ and 
$\sq\le j_r(P)\le \sq+1$ at $P\not\in X(\fq)$. In this case we have 
$j_r(P)=\sq$: in fact, if $j_r(P)=\sq+1$, by (\ref{last}), we would have 
$(\sq+1)P\sim \sq P+\frx(P)$, i.e. $P\sim \frx(P)$ and so $g=0$, a 
contradiction. Therefore $\nu_{r-1}=\epsilon_r=\sq$. Now, the invariants of $\cD$ 
also 
satisfy properties (\ref{eq2}) and (\ref{eq5}) above.  The case $i=r-1$ of   
(\ref{eq5}) gives $j_1(P)=1$ for each $P\in X(\fq)$. Hence by (\ref{eq2}) we find
that 
$m_{r-1}(P)=\sq$ for each $P\in X(\fq)$. It follows that $\cD$ is simple. 
Now by using the hypothesis on $g$ and 
Castelnuovo's genus bound we find that $r=2$ (cf. \cite[Claim 1]{ft}).

Then, by using the ramification divisor of $\cD$, we have that
$g=\sq(\sq-1)/2$ (see \cite[\S3]{ft}). Fix $P_0\in X(\fq)$,  
$x, y\in \fq(X)$ with $\dinf(x)=\sq P_0$, and $\dinf(y)=(1+\sq)P_0$. Then
$H(P_0)=\langle \sq,\sq+1\rangle$ and consequently $\df(dx)=(2g-2)P_0 (*)$, 
because $H(P_0)$ is symmetric. 

By $\nu_1=\sq$ we have an equation of type 
\begin{equation}\label{last1}
y^q-y=f(x^q-x),
\end{equation}
with $f:=D^{1}y$ (derivation with respect to $x$), and by $\epsilon_2=\sq$ we have 
the following rank 2 matrices
$$
\left( \begin{array}{ccc}
1 & x   & y\\
0 & 1   & D^{(1)}y\\
0 & 0   & D^{(j)}y
\end{array} \right),
\quad 2\le j<\epsilon_2=\sq.
$$
By $(*)$ and (\ref{last1}), $\dinf(f)=qP_0$.  
Then $D^{(j)}y=0$ for $2\le j<\sq$ and from (\ref{last1}) we
have that $D^{(j)}f=0$ for $1\le j<\sq$. Then by 
\cite[Satz 10]{hasse}, there exists $f_1\in \fq(X)$ such that  
$f=f_1^{\sq}$. Therefore $f_1=ax+b$ with $a, b\in
\fq$, $a\neq0$, and after some $\fq$-linear transformations we obtain an
equation of type
$$
y_1^{\sq}+y_1-x_1^{\sq+1}=(y_1^{\sq}+y_1-x_1^{\sq+1})^{\sq},
$$
with $x_1, y_1\in \fq(X)$. Now the proof follows.
\begin{remark*} Similar computations to those of the Remarks in \S3 can be
obtained in this case. In particular, using the divisor $\cD$, one can
show that Hermitian curves are non-classical curves 
(see also \cite{gv}). We mention that the 
first examples of non-classical curves were obtained from some Hermitian
curves (see \cite{sch}).
\end{remark*}
\smallskip

\begin{center}
{\bf Appendix:} On the Suzuki-Tits ovoid
\end{center}
\smallskip

For $s\in \N$, let $q_0:=2^s$ and $q:=2q_0$. It is well known that the
Suzuki-Tits ovoid $\cO$ can be represented in $\P^4(\fq)$ as
$$
\cO=\{(1:a:b:f(a,b):af(a,b)+b^2): a, b \in \fq\}\cup\{(0:0:0:0:0:1)\},
$$
where $f(a,b):=a^{2q_0+1}+b^{2q_0}$ (see \cite{tits}, \cite[p.3]{pent}) 

Let $X$ be the Deligne-Lusztig curve associated to $Sz(q)$ and 
$\cD=|(q+2q_0+1)P_0|$, $P_0\in X(\fq)$ (see \S1). By Remark \S3(3),
we can associate to $\cD$ a morphism $\pi=(1:x:y:z:w)$ whose coordinates
satisfy $\dinf(x)=qP_0$, $\dinf(y)=(q+q_0)P_0$, $\dinf(z)=(q+2q_0)P_0$ and
$\dinf(w)=q+2q_0+1$. 
\begin{claim*} (A. Cossidente)\quad 
$\cO=\pi(X(\fq))$.
\end{claim*}
\begin{proof} We have $\pi(P_0)=(0:0:0:0:1)$; we can choose 
$x$ and $y$ satisfying the plane equation (given in \S0) of the curve, 
$z:= x^{2q_0+1}+y^{2q_0}$, and $w:=
xy^{2q_0}+z^{2q_0}=xy^{2q_0}+x^{2q+2q_0}+y^{2q}$
(cf. \cite[S1.7]{hsti}). For $P\in X(\fq)\setminus\{P_0\}$ set $a:=x(P)$, 
$b:=y(P)$, and $f(a,b):= z(a,b)$. Then $w(a,b)=af(a,b)+b^2$ and we are
done.
\end{proof}
\begin{remark*} The morphism $\pi$ is an embedding. Indeed, since
$j_1(P)=1$ for each $P$ (cf. Remarks \S3(2)(4)), it is enough to see that
$\pi$ is injective. By (\ref{eq1}), the points $P$ where $\pi$ could not
be injective satisfy: $P\not\in X(\fq)$$, \frx^3(P)=P$ or $\frx^2(P)=P$.
Now from the Zeta function of $X$ one sees that $\#X(\mathbb
F_{q^3})=\#X(\mathbb F_{q^2})=\#X(\fq)$, and the remark follows.
\end{remark*}

\end{document}